\documentclass[12pt,a4paper]{amsart}

\usepackage[utf8]{inputenc}
\usepackage[T1]{fontenc}
\usepackage[english]{babel}
\usepackage{cite}

\usepackage{indentfirst}
\usepackage{amssymb}
\usepackage{amsfonts}
\usepackage{amsmath}
\usepackage{amsthm}
\usepackage[top=2.5cm, bottom=2.5cm, left=2.5cm, right=2.5cm]{geometry}
\usepackage{amsopn}
\usepackage[colorlinks]{hyperref}

\theoremstyle{plain}
\newtheorem{thm}{Theorem}

\newtheorem{lem}[thm]{Lemma}
\newtheorem{prop}[thm]{Proposition}

\theoremstyle{definition}
\newtheorem{defn}[thm]{Definition}
\newtheorem{example}[thm]{Example}
\newtheorem*{example*}{Example}

\newtheorem*{rem*}{Remark}

\newcommand{\R}{\mathbb{R}}

\newcommand{\udimA}{\overline{\text{dim}}_{A}}

\DeclareMathOperator{\diam}{diam}

\DeclareMathOperator{\ucodim}{\underline{co\,dim}}
\DeclareMathOperator{\ocodim}{\overline{co\,dim}}

\DeclareMathOperator{\dist}{dist}

\DeclareMathOperator{\WLSC}{WLSC}
\DeclareMathOperator{\WUSC}{WUSC}

\DeclareSymbolFont{bbsymbol}{U}{bbold}{m}{n}
\DeclareMathSymbol{\ind}{\mathbin}{bbsymbol}{'061}

\title[On density of compactly supported functions in fractional Sobolev spaces]{On density of compactly supported smooth functions in fractional Sobolev spaces}
\author[B{.} Dyda]{Bart{\l}omiej Dyda}
\author[M{.} Kijaczko]{Micha\l{} Kijaczko}


\address[B.D. and M.K.]{Faculty of Pure and Applied Mathematics\\ Wroc{\l}aw University 
	of Science and Technology\\
	Wybrze\.ze Wyspia\'nskiego 27,
	50-370 Wroc{\l}aw, Poland
}
\email{bdyda@pwr.edu.pl\qquad dyda@math.uni-bielefeld.de}
\email{michal.kijaczko@pwr.edu.pl}
\thanks{B.D. was partially supported by grant NCN 2015/18/E/ST1/00239.}

\begin{document}

\begin{abstract}
We describe some sufficient conditions, under which smooth and compactly supported functions are or are not dense in the fractional Sobolev space $W^{s,p}(\Omega)$ for an open, bounded set $\Omega\subset\mathbb{R}^{d}$. The density property is closely related to the lower and upper Assouad codimension of the boundary of $\Omega$. We also describe explicitly the closure of $C_{c}^{\infty}(\Omega)$ in $W^{s,p}(\Omega)$ under some mild assumptions about the geometry of $\Omega$. Finally, we prove a~variant of a~fractional order Hardy inequality.
\end{abstract}

	\maketitle
	\tableofcontents 
	
	\section{Introduction}
	We discuss the problem of density of compactly supported smooth functions in the fractional Sobolev space $W^{s,p}(\Omega)$, which is well known to hold when $\Omega$ is a bounded Lipschitz domain and $sp\leq 1$ \cite[Theorem 1.4.2.4]{MR775683},\cite[Theorem 3.4.3]{MR3024598}.
        We extend this result to bounded, plump open sets with a~dimension of the boundary satisfying certain inequalities. To this end, we use the Assouad dimensions and codimensions.\
        We  also describe explicitly the closure of $C_{c}^{\infty}(\Omega)$ in the fractional Sobolev space, provided that $\Omega$ satisfies the fractional Hardy inequality.

	Let $\Omega\subset\mathbb{R}^{d}$ be an open set. Let $0<s<1$ and $1\leq p<\infty$. We recall that the \emph{fractional Sobolev space} is defined as 
	$$
	W^{s,p}(\Omega)=\left\{f\in L^{p}(\Omega):\int_{\Omega}\int_{\Omega}\frac{|f(x)-f(y)|^{p}}{|x-y|^{d+sp}}\,dy\,dx<\infty\right\}.
	$$
	This is a Banach space endowed with the norm 
	$$
	\|f\|_{W^{s,p}(\Omega)}=\|f\|_{L^{p}(\Omega)}+[f]_{W^{s,p}(\Omega)},
	$$
	where $[f]_{W^{s,p}(\Omega)}=\left(\int_{\Omega}\int_{\Omega}\frac{|f(x)-f(y)|^{p}}{|x-y|^{d+sp}}\,dy\,dx\right)^{1/p}$ is called the \emph{Gagliardo seminorm}.
        Throughout the paper we consider only real-valued functions, but we note that all results are clearly valid also for complex-valued functions, by means of decomposing them into a~sum of real and imaginary part.
	\begin{defn}\label{defn1}
	By $W_{0}^{s,p}(\Omega)$ we denote the closure of $C_{c}^{\infty}(\Omega)$ (the space of all smooth functions with compact support in~$\Omega$) in $W^{s,p}(\Omega)$ with respect to the Sobolev norm.
	\end{defn}

        The following theorem is our main result on the connection between $W_{0}^{s,p}(\Omega)$
        and $W^{s,p}(\Omega)$. For the relevant geometric definitions, we refer the Reader to Section~\ref{sec:geom}. Here we only note that for bounded Lipschitz domains one has
        $\ucodim_{A}(\partial\Omega) = \ocodim_{A}(\partial\Omega) = 1$ and the other geometrical assumptions of Theorem~\ref{tw1} do hold (that is, bounded Lipschitz domains are $(d-1)$-homogeneous and $\kappa$-plump), hence the classical case is included.
\begin{thm}\label{tw1}
Let $\Omega\subset \R^d$ be a nonempty bounded open set, let $0<s<1$ and $1\leq p <\infty$.\\
(I) If  $sp<\ucodim_{A}(\partial\Omega)$,
then $W_{0}^{s,p}(\Omega)=W^{s,p}(\Omega)$.\\
(II)  If $\Omega$ is a $(d-sp)$-homogeneous set,
$sp=\ucodim_{A}(\partial\Omega)$
and $p>1$, 
then $W_{0}^{s,p}(\Omega)=W^{s,p}(\Omega)$.\\
(III) If $\Omega$ is $\kappa$-plump and $sp>\ocodim_{A}(\partial\Omega)$,
then $W_{0}^{s,p}(\Omega)\neq W^{s,p}(\Omega)$.\\
\end{thm}
We remark that a result similar to the part (I) and (III) in the Theorem \ref{tw1} was obtained by Caetano in \cite{MR1792288} in the context of Besov spaces and Triebel--Lizorkin spaces, but with the Minkowski dimension instead of Assouad dimension. That result is not directly comparable with ours, as for less regular domains spaces $W^{s,p}$ do not necessiraly coincide with the appropriate Triebel--Lizorkin spaces.
We  refer the Reader to \cite{MR3987216} for a~discussion on the space $W^{s,p}_0$ and different similarly defined spaces. We also want to mention that analogous, but slightly different problems were considered in \cite{MR3310082} (spaces of functions vanishing outside $\Omega$), \cite{MR3420496} (the weighted case) and \cite{MR3989177} (spaces with variable exponents).

In the case (III) above we also obtain the following characterisation of the space $W_{0}^{s,p}(\Omega)$. For the proof, see Section~\ref{sec:W0}.

\begin{thm}\label{twW0}
  Let $0<s<1$ and $1\leq p < \infty$.
  Suppose that $\Omega\neq\emptyset$ is a bounded, open $\kappa$-plump set.
  If  $\ocodim_{A}(\partial\Omega)<sp$, then
  \begin{equation}\label{eq:W0isW}
  W_{0}^{s,p}(\Omega) = \left\{f\in W^{s,p}(\Omega):\int_{\Omega}\frac{|f(x)|^{p}}{\dist(x,\partial \Omega)^{sp}}\,dx<\infty\right\}.
  \end{equation}
\end{thm}

  In the case (I) of Theorem~\ref{tw1} equality \eqref{eq:W0isW} also holds,
  or in other words, we have an inclusion between the Sobolev and weighted $L^p$ space,
  $W^{s,p}(\Omega) \subset L^p(\Omega, \dist(x,\partial \Omega)^{-sp})$.
 This fact is made quantitative in the next theorem; for its proof, see Section~\ref{sec:W0} as well.
  \begin{thm}\label{weakH}
     Let $0<s<1$ and $1\leq p < \infty$.
  Suppose that $\Omega\neq\emptyset$ is a bounded, open $\kappa$-plump set.
  If  $\ucodim_{A}(\partial\Omega)>sp$, then there exists
a~constant $c$ such that
\begin{equation}\label{eq:weakH}
  \int_{\Omega}\frac{|f(x)|^{p}}{\dist(x,\partial \Omega)^{sp}}\,dx
\leq c \|f\|_{W^{s,p}(\Omega)}^p < \infty, \quad \text{for all $f\in W^{s,p}(\Omega)$}.
\end{equation}
\end{thm}

  Theorem~\ref{twW0} and \ref{weakH} have classical (non-fractional) counterparts,
  see \cite[Example 9.11]{MR664599} or \cite{MR1470421}.

Finally, we extend the results of \cite[Theorem 1, Corollary 3]{MR3237044}.
Namely, we  prove the case (T') in the following version of the fractional Hardy inequality.
For the definitions of the conditions
$\WLSC$ and $\WUSC$ we refer the reader to the Appendix,
while  the plumpness and Assouad dimensions are defined in Section~\ref{sec:geom}.
We would also like to note that a~special case of (T') (assuming in particular $p=2$) was proved in \cite[Lemma 3.32]{MR1742312} and \cite{MR1993864}.

\begin{thm}[\cite{MR3237044} in cases (T) and (F)]\label{thm.dv}
Let $0<p<\infty$, $H\in (0,1]$ and $\eta\in\R$.
Suppose $\Omega\not=\emptyset$ is a proper $\kappa$-plump open set in $\R^d$
and $\phi:(0,\infty)\to(0,\infty)$ is a function
 so that either condition (T), or condition (T'), or condition (F) holds
\begin{itemize}
\item[(T)]  $\eta + \overline{\mathrm{dim}}_A(\partial \Omega) - d<0$, $\Omega$ is unbounded,
$\phi\in \WUSC(\eta,0,H^{-1})$,
\item[(T')]  $\eta + \overline{\mathrm{dim}}_A(\partial \Omega) - d<0$, $\Omega$ is bounded,
$\phi\in \WUSC(\eta,0,H^{-1})$,
\item[(F)]  $\eta + \underline{\mathrm{dim}}_A(\partial \Omega) - d > 0$, $\Omega$ is bounded or $\partial \Omega$ is unbounded, and  $\phi\in \WLSC(\eta,0,H)$.
\end{itemize}
Then there exist constants $c=c(d,s,p,\Omega,\phi)$ and $R$ such that the following inequality
\begin{equation}\label{eq:fhi}
 \int_\Omega \frac{|u(x)|^p}{\phi(d_{\Omega}(x))}\,dx \leq
c \int_\Omega\!\int_{\Omega\cap B(x,Rd_{\Omega}(x))}
 \frac{|u(x)-u(y)|^p}{\phi(d_{\Omega}(x))d_{\Omega}(x)^d}\,dy\,dx\, + c\xi \|u\|_{L^p(\Omega)}^p,
\end{equation}
holds for all measurable functions $u$ for which the left hand side is finite, with $\xi=0$ in the cases $(T)$ and $(F)$ and $\xi=1$ in the case $(T')$.
\end{thm}
There is a huge literature about fractional Hardy inequalities; we refer the Reader to \cite{MR2085428, MR3237044, MR3277052} and the references therein.
We would also like to draw Reader's attention to a paper \cite{MR1810244} from 1999 by Farman Mamedov.
This not very well-known paper is one of the first to deal with multidimensional fractional order Hardy inequalities.

The authors would like to thank Lorenzo Brasco for helpful discussions on the subject, in particular for providing a part of the proof of Theorem~\ref{tw1}, and the anonymous referee for numerous comments which led to an
improvement of the manuscript.

\section{Geometrical definitions}\label{sec:geom}
We denote  the distance from $x\in \R^d$ to a~set $E\subset \R^d$ by
$\dist(x,E)=\displaystyle\inf_{y\in E}|x-y|$;
for open sets $\Omega\subset \R^d$ we write
 $d_{\Omega}(x)=\text{dist}(x,\partial\Omega)$.

\begin{defn}\label{defn3}
  Let $r>0$. For open sets $\Omega\subset \R^d$
  we define the \emph{inner tubular neighbourhood}  of $\Omega$ as
\[
\Omega_{r}=\left\{x\in\Omega:d_{\Omega}(x)\leq r\right\},
\]
and for arbitrary sets $E\subset \R^d$ we define the \emph{tubular neighbourhood} of $E$ as
$$
\widetilde{E}_{r}=\left\{x\in\mathbb{R}^{d}: \dist(x, E) \leq r\right\}.
$$

\end{defn}
\begin{defn}
\cite[Section 3]{MR3205534} Let $E\subset\mathbb{R}^{d}$. The \emph{lower Assouad codimension} $\underline{\text{co\,dim}}_{A}(E)$ is defined as the supremum of all $q\geq 0$, for which there exists a constant $C=C(q)\geq 1$ such that for all $x\in E$ and $0<r<R<\diam E$ it holds
$$
\left|\widetilde{E}_{r}\cap B(x,R)\right|\leq C\left|B(x,R)\right|\left(\frac{r}{R}\right)^{q}.
$$
Conversely, the \emph{upper Assouad codimension} $\overline{\text{co\,dim}}_{A}(E)$ is defined as the infimum of all $s\geq 0$, for which there exists a constant $c=c(s)>0$ such that for all $x\in E$ and $0<r<R<\diam E$ it holds
$$
\left|\widetilde{E}_{r}\cap B(x,R)\right|\geq c\left|B(x,R)\right|\left(\frac{r}{R}\right)^{s}.
$$
\end{defn}

We remark that having strict inequality $R<\diam E$ above
makes the definitions applicable also for unbounded sets $E$;
for bounded sets $E$ we could have $R\leq \diam E$.

In Euclidean space $\mathbb{R}^{d}$ we have $\underline{\text{dim}}_{A}(E)=d-\overline{\text{co\,dim}}_{A}(E)$,  $\overline{\text{dim}}_{A}(E)=d-\underline{\text{co\,dim}}_{A}(E)$, where $\underline{\text{dim}}_{A}(E) $ and $\overline{\text{dim}}_{A}(E)$ denote respectively the well known lower and upper Assouad dimension -- see for example \cite[Section 2]{MR3205534} for this result. Recall that the upper Assouad dimension of a given set $E$ is defined as the infimum of all exponents $s\geq 0$ for which there exists a constant $C=C(s)\geq 1$ such that for all $x\in E$ and $0<r<R<\diam E$ the ball $B(x,R)\cap E$ can be covered by at most $C(R/r)^s$ balls with radius $r$, centered at $E$. Analogously, the lower Assouad dimension is characterized by the supremum of all exponents $t\geq 0$ for which there is a constant $c=c(t)>0$ such that the ball $B(x,R)\cap E$ can be covered by at least c$(R/r)^t$ balls with radius $r$ and centered at $E$. If $\underline{\text{co\,dim}}_{A}(E)=\overline{\text{co\,dim}}_{A}(E),$ we simply denote it by $\text{co\,dim}_{A}(E)$.

	We recall a geometric notion from \cite{10.2748/tmj/1178228081}.
	\begin{defn}
		A set $E\subset \R^d$ is {\em $\kappa$-plump}
		with $\kappa\in (0,1)$ if, for each $0<r< \diam(E)$ and each $x\in \overline{E}$, there
		is $z\in \overline{B}(x,r)$ such that
		$B(z,\kappa r)\subset E$.\end{defn}

        Following \cite[Theorem A.12]{MR1608518}, we define a notion of $\sigma$-homogenity.
\begin{defn}
Let $E\subset\R^d$ and let $V(E,x,\lambda,r)=\{y\in\R^d:\dist(y,E)\leq r, |x-y|\leq \lambda r\}$. We say that E is $\sigma$-\emph{homogeneous}, if there exists a~constant $L$ such that
$$
|V(E,x,\lambda,r)|\leq Lr^{d}\lambda^{\sigma}
$$
for all $x\in E$, $\lambda \geq 1$ and $r>0$.
\end{defn}
If $0<r<R<\diam(E)$, then taking $\lambda=R/r$ in the definition gives
\[
  \left|\widetilde{E}_{r} \cap B(x,R)\right| = \left|V\Big(E,x,\frac{R}{r},r\Big)\right| \leq C\left|B(x,R)\right|\left(\frac{r}{R}\right)^{d-\sigma},
\]
where $C=C(d,E)$ is a constant. This means that if $\ucodim_{A}(E)=s$, then  $(d-s)$-homogeneous sets are precisely these sets $E$, for which the supremum in the definition of the lower Assouad codimension is attained. For the definition of the concept of homogenity 
from a different point of view the Reader may also see \cite[Definition 3.2]{MR1608518}.

Finally, let us note that for example in part I of Theorem~\ref{tw1}, we need the assumption  $sp<\ucodim_{A}(\partial\Omega)$ only to obtain the bound \eqref{eq:tubu}.
For that a~slightly weaker assumption in terms of  Minkowski (co)dimension would suffice, however, we need Assouad (co)dimensions for other parts of the paper and therefore we prefer to use only them. Let us only recall that the upper Minkowski dimension of a set $E\subset\R^d$ is defined as
$$
\overline{\text{dim}}_{M}(E)=\inf\{s\geq0: \limsup_{r\rightarrow0}\left|\widetilde{E}_{r}\right|r^{d-s}=0\},
$$
see for example \cite[Section 2]{MR4144553}. The statement of the part (I) of Theorem \ref{tw1} remains true if we assume that $sp<d-\overline{\text{dim}}_{M}(\partial\Omega)$.
\section{Lemmas}

The following lemma is the key to our further computations. We recall that $\Omega_{\frac{3}{n}}$
appearing in \eqref{eq:fvn} is the inner tubular neighbourhood of $\Omega$, see Definition~\ref{defn3}.
\begin{lem}\label{lem1}
  Let
\[
  v_{n}(x)= \max\left\{\min\left\{2 - nd_{\Omega}(x),1\right\},0\right\}=\left\{ \begin{array}{ll}
1 & \textrm{when $d_{\Omega}(x)\leq 1/n$,}\\
2 - nd_{\Omega}(x) & \textrm{when $ 1/n < d_{\Omega}(x)\leq2/n$,}\\
0 & \textrm{when $d_{\Omega}(x)>2/n$.}
  \end{array} \right.
  \]
  There exists a constant $C=C(d,s,p,\Omega)>0$ such that the following inequality holds
  for all functions $f\in W^{s,p}(\Omega)$
  \begin{equation}\label{eq:fvn}
    [fv_{n}]^{p}_{W^{s,p}(\Omega)}
    \leq Cn^{sp}\int_{\Omega_{\frac{3}{n}}}|f(x)|^{p}\,dx+C\int_{\Omega_{\frac{3}{n}}}\int_{\Omega_{\frac{3}{n}}}\frac{|f(x)-f(y)|^{p}}{|x-y|^{d+sp}}\,dy\,dx.
\end{equation}
\end{lem}
\begin{proof}
  Fix $f\in W^{s,p}(\Omega)$ and define $f_{n}=fv_{n}$. We have
\begin{align*}
[f_{n}]^{p}_{W^{s,p}(\Omega)}&=\int_{\Omega}\int_{\Omega}\frac{|f(x)v_{n}(x)-f(y)v_{n}(y)|^{p}}{|x-y|^{d+sp}}\,dy\,dx\\
&=\int_{\Omega_{\frac{3}{n}}}\int_{\Omega_{\frac{3}{n}}}\frac{|f(x)v_{n}(x)-f(y)v_{n}(y)|^{p}}{|x-y|^{d+sp}}\,dy\,dx\\
&+2\int_{\Omega_{\frac{2}{n}}}\int_{\Omega\setminus\Omega_{\frac{3}{n}}}\frac{|f(x)v_{n}(x)|^{p}}{|x-y|^{d+sp}}\,dy\,dx\\
&=:J_{1}+2J_{2}.
\end{align*}
First we estimate $J_{1}$,
\begin{align*}
2^{1-p}J_{1}&\leq\int_{\Omega_{\frac{3}{n}}}\int_{\Omega_{\frac{3}{n}}}\frac{|f(x)|^{p}|v_{n}(x)-v_{n}(y)|^{p}}{|x-y|^{d+sp}}\,dy\,dx\\
&+\int_{\Omega_{\frac{3}{n}}}\int_{\Omega_{\frac{3}{n}}}\frac{|v_{n}(y)|^{p}|f(x)-f(y)|^{p}}{|x-y|^{d+sp}}\,dy\,dx\\
&=:K_{1}+K_{2}.
\end{align*}
Since $|v_n|\leq 1$, we obtain
$$
K_{2}\leq \int_{\Omega_{\frac{3}{n}}}\int_{\Omega_{\frac{3}{n}}}\frac{|f(x)-f(y)|^{p}}{|x-y|^{d+sp}}\,dy\,dx.
$$
Furthermore, $|v_{n}(x)-v_{n}(y)|\leq \min\{1,n|x-y|\}$, hence, for $K_{1}$ we can compute that
\begin{align*}
K_{1}&\leq \int_{\Omega_{\frac{3}{n}}}\int_{\Omega_{\frac{3}{n}}}\frac{|f(x)|^{p}\left(\min\{1,n|x-y|\}\right)^{p}}{|x-y|^{d+sp}}\,dy\,dx\\  
&\leq\int_{\Omega_{\frac{3}{n}}}|f(x)|^{p}\,dx\int_{|x-y|>1/n}\frac{\,dy}{|x-y|^{d+sp}}+n^{p}\int_{\Omega_{\frac{3}{n}}}|f(x)|^{p}\,dx\int_{|x-y|<1/n}\frac{\,dy}{|x-y|^{d-(1-s)p}}\\
&\leq Cn^{sp}\int_{\Omega_{\frac{3}{n}}}|f(x)|^{p}\,dx.
\end{align*}
Since $|v_{n}|\leq 1$, for $J_{2}$ we have 
\begin{align*}
J_{2}=&\int_{\Omega_{\frac{2}{n}}}\int_{\Omega\setminus\Omega_{\frac{3}{n}}}\frac{|f(x)|^{p}|v_{n}(x)|^{p}}{|x-y|^{d+sp}}\,dy\,dx\\
&\leq \int_{\Omega_{\frac{2}{n}}}|f(x)|^{p}\,dx\int_{\Omega\setminus\Omega_{\frac{3}{n}}}\frac{\,dy}{|x-y|^{d+sp}}\\
&\leq \int_{\Omega_{\frac{2}{n}}}|f(x)|^{p}\,dx\int_{B(x,1/n)^{c}}\frac{\,dy}{|x-y|^{d+sp}}\\
&\leq Cn^{sp}\int_{\Omega_{\frac{2}{n}}}|f(x)|^{p}\,dx.
\end{align*}
Hence, we obtain for some (new) constant $C$ that 
\begin{align*}
[f_{n}]^{p}_{W^{s,p}(\Omega)}&\leq Cn^{sp}\int_{\Omega_{\frac{3}{n}}}|f(x)|^{p}\,dx+C\int_{\Omega_{\frac{3}{n}}}\int_{\Omega_{\frac{3}{n}}}\frac{|f(x)-f(y)|^{p}}{|x-y|^{d+sp}}\,dy\,dx.\qedhere
\end{align*}
\end{proof}

\begin{defn}\label{defn2}
  By $W_{c}^{s,p}(\Omega)$ we denote the closure of all compactly supported functions in $W^{s,p}(\Omega)$ (not necessarily smooth) with respect to the Sobolev norm.
\end{defn}
The key property, which allows us to get rid of the smoothness and rely only on the compactness of the support, is the result below.
\begin{prop}\label{prop1}
  We have $W_{0}^{s,p}(\Omega)=W_{c}^{s,p}(\Omega)$.
\end{prop}
\begin{proof}
  This is a straightforward consequence of \cite[Proposition 2 and proof of Theorem 8]{MR4190640}.
\end{proof}
It turns out that to prove the density of compactly supported functions in the fractional Sobolev space, we only need to find a sequence which approximates the function $\ind_{\Omega}$ (the indicator of $\Omega$).

\begin{lem}\label{lem2}
Let $\Omega$ be an open set such that $|\Omega|<\infty$. We have
$$
W_{0}^{s,p}(\Omega)=W^{s,p}(\Omega) \Longleftrightarrow \ind_{\Omega}\in W_{0}^{s,p}(\Omega)
$$
\end{lem}
\begin{proof}
  Implication ``$\implies$'' is obvious, therefore we proceed to prove the implication from right to left.
According to Proposition \ref{prop1}, we need to prove that if the function $\ind_{\Omega}$ can be approximated by some family of functions $g_{n}\in W_{c}^{s,p}(\Omega)$, then every function $f\in W^{s,p}(\Omega)$ can be approximated by functions from $W_{c}^{s,p}(\Omega)$. Since $L^{\infty}(\Omega)\cap W^{s,p}(\Omega)$ is dense in $W^{s,p}(\Omega)$ (because the truncated functions $f^{N}=\min\left\{\max\left\{f,-N\right\},N\right\}$ tend to $f$ in $W^{s,p}(\Omega)$, as $N\longrightarrow\infty$), we may assume that $f\in L^{\infty}(\Omega)$. Moreover, we may also assume that $0\leq g_{n}\leq 1$, because if $g_{n}\longrightarrow \ind_{\Omega}$ in $W^{s,p}(\Omega)$, then also $\widetilde{g}_{n}=\max\{\min\{g_{n},1\},0\}\longrightarrow\ind_{\Omega}$, since we have $|\widetilde{g}_{n}(x)-\widetilde{g}_{n}(y)|\leq |g_{n}(x)-g_{n}(y)|$.

Define $f_{n}=fg_{n}\in W_{c}^{s,p}(\Omega)$. Observe that
\begin{align*}
[f-f_{n}]^{p}_{W^{s,p}(\Omega)}
&=\int_{\Omega}\int_{\Omega}\frac{|f(x)(1-g_{n}(x))-f(y)(1-g_{n}(y))|^{p}}{|x-y|^{d+sp}}\,dy\,dx\\
&\leq 2^{p-1}\int_{\Omega}\int_{\Omega}\frac{|f(x)|^{p}|g_{n}(x)-g_{n}(y)|^{p}}{|x-y|^{d+sp}}\,dy\,dx\\
&+2^{p-1}\int_{\Omega}\int_{\Omega}\frac{|1-g_{n}(y)|^{p}|f(x)-f(y)|^{p}}{|x-y|^{d+sp}}\,dy\,dx\\
&\leq 2^{p-1}\|f\|^{p}_{\infty}[g_{n}]^{p}_{W^{s,p}(\Omega)}\\
&+2^{p-1}\int_{\Omega}\int_{\Omega}\frac{|1-g_{n}(y)|^{p}|f(x)-f(y)|^{p}}{|x-y|^{d+sp}}\,dy\,dx.
\end{align*}
Since $g_{n}\longrightarrow \ind_{\Omega}$ in $L^{p}(\Omega)$, there is a subsequence $g_{n_{k}}\longrightarrow \ind_{\Omega}$ almost everywhere. Hence, for such a subsequence we have
\begin{align*}
[f-f_{n_{k}}]^{p}_{W^{s,p}(\Omega)}&\leq 2^{p-1}\|f\|^{p}_{\infty}[g_{n_{k}}]^{p}_{W^{s,p}(\Omega)}\\
&+2^{p-1}\int_{\Omega}\int_{\Omega}\frac{|1-g_{n_{k}}(y)|^{p}|f(x)-f(y)|^{p}}{|x-y|^{d+sp}}\,dy\,dx.
\end{align*}
The first term above is convergent to $0$, since $g_{n_{k}}\longrightarrow\ind_{\Omega}$ in $W^{s,p}(\Omega)$. The convergence of the second term follows from Lebesgue dominated convergence theorem. Moreover, it is trivial to show that $f_{n}\longrightarrow f$ in $L^{p}(\Omega)$ and hence the proof is finished.
\end{proof}

\section{Proof of Theorem~\ref{tw1}}

\begin{proof}[Proof of Theorem~\ref{tw1}, case I]
  According to Lemma \ref{lem2}, we only need to prove that the function $f=\ind_{\Omega}$ can be approximated by compactly supported functions. Let $f_{n}=fv_{n}$, where $v_{n}$ is as in the Lemma \ref{lem1} and let $\underline{d}=\ucodim_{A}(\partial\Omega)$. By Lemma \ref{lem1} (note that
in this case the second term in inequality \eqref{eq:fvn} is $0$) we have
\begin{align*}
[f_{n}]^{p}_{W^{s,p}(\Omega}\leq Cn^{sp}\int_{\Omega_{\frac{3}{n}}}\,dx=Cn^{sp}\left|\Omega_{\frac{3}{n}}\right|.
\end{align*}
If $sp<\underline{d}$, then, by the definition of lower Assouad codimension, for every $\varepsilon>0$ we have
\begin{equation}\label{eq:tubu}
  \left|\Omega_{\frac{3}{n}}\right|\leq C'\left(\frac{1}{n}\right)^{\underline{d}-\varepsilon}.
\end{equation}
Hence, for some new constant C we have
\begin{align*}
    [f_{n}]^{p}_{W^{s,p}(\Omega)}\leq Cn^{sp}n^{\varepsilon-\underline{d}}\longrightarrow 0,
\end{align*}
when $n\longrightarrow\infty$, by choosing $0<\varepsilon<\underline{d}-sp$, which is feasible thanks to our assumption.
\end{proof}

\begin{proof}[Proof of Theorem~\ref{tw1}, case II]
  We proceed like in the above proof of the first part of the Theorem \ref{tw1} and obtain
  \[
    [f_{n}]^{p}_{W^{s,p}(\Omega)}\leq Cn^{sp}\left|\Omega_{\frac{3}{n}}\right|.
    \]
    Since $\Omega$ is $(d-sp)$-homogeneous and $\ucodim_{A}(\partial\Omega)=sp$, then it follows that $\left|\Omega_{\frac{3}{n}}\right|\leq C'n^{-sp}$ and, in consequence, the sequence $\{f_{n}\}_{n\in\mathbb{N}}$ is bounded in $W^{s,p}(\Omega)$.

    The following argument was kindly pointed out to us by Lorenzo Brasco,
    see also \cite[Theorem 4.4]{MR4225499} for a~similar argument.
    It is well known that for $p>1$ the space $W^{s,p}(\Omega)$ is reflexive. Hence, by Banach--Alaoglu and Eberlein--\v{S}mulian theorem, there exists a subsequence $\{f_{n_{k}}\}_{k\in\mathbb{N}}$ weakly convergent to some $f$. Since $W_{0}^{s,p}(\Omega)$ is both closed and convex subset of $W^{s,p}(\Omega)$, by \cite[Theorem 2.3.6]{MR3309383} it is also weakly closed, so we have $f\in W_{0}^{s,p}(\Omega)$. Then it suffices to see that  $f=\ind_{\Omega}$ by the uniqueness of the limit, since $f_{n_{k}}$ strongly converges to $\ind_{\Omega}$ in $L^{p}(\Omega)$. This ends the proof.
\end{proof}

\begin{proof}[Proof of Theorem~\ref{tw1}, case III]
 Let $\overline{d}=\overline{\text{co\,dim}}_{A}(\partial\Omega)$. We will show that the indicator of $\Omega$ cannot be approximated by functions with compact support. Indeed, let $u_{n}$ be any sequence of compactly supported function such that $\|u_{n}-\ind_{\Omega}\|_{W^{s,p}(\Omega)}\longrightarrow 0.$ In particular $u_{n}\longrightarrow \ind_{\Omega}$ in $L^{p}(\Omega),$ so there is a subsequence $u_{n_{k}}$ convergent almost everywhere to $\ind_{\Omega}$. If $sp>\overline{d}$, we can use the fractional Hardy inequality from \cite[Corollary 3]{MR3237044} in the case (F) with $\beta=0$ to obtain

\begin{align*}
[u_{n_{k}}-\ind_{\Omega}]^{p}_{W^{s,p}(\Omega)}&=[u_{n_{k}}]^{p}_{W^{s,p}(\Omega)}=\int_{\Omega}\int_{\Omega}\frac{|u_{n_{k}}(x)-u_{n_{k}}(y)|^{p}}{|x-y|^{d+sp}}\,dy\,dx\\
&\geq c\int_{\Omega}\frac{|u_{n_{k}}(x)|^{p}}{d_{\Omega}(x)^{sp}}\,dx.
\end{align*}
By Fatou's lemma,
\begin{align*}
0=\lim_{k\rightarrow\infty}[u_{n_{k}}]^{p}_{W^{s,p}(\Omega)}&\geq c\int_{\Omega}\liminf_{k\rightarrow\infty}\frac{|u_{n_{k}}(x)|^{p}}{d_{\Omega}(x)^{sp}}\,dx\\
&=c\int_{\Omega}\frac{\,dx}{d_{\Omega}(x)^{sp}}>0.
\end{align*}
We obtain a contradiction.
\end{proof}

\begin{example}(Lipschitz domains)
Let $\Omega$ be a bounded Lipschitz domain. In this case  we have $\text{co\,dim}_{A}(\partial\Omega)=1$ and, by the cone property, $\left|\Omega_{r}\right|=O(r)$, hence, Theorem~\ref{tw1} generalises the classical result \cite[Theorem 1.4.2.4]{MR775683}.
\end{example}
\begin{example}(Koch snowflake)
  Let $\Omega\subset\mathbb{R}^{2}$ denote the domain bounded by the Koch snowflake. It is well known that the Hausdorff dimension of the Koch \emph{curve} is $\frac{\log4}{\log 3}$.
  Thus also its Assouad dimension is $\frac{\log4}{\log 3}$, since it is a~self-similar set
  satisfying open set condition, see \cite[Corollary~2.11]{MR3267023}.
  The Koch snowflake is a~finite union of copies of Koch curves, therefore its Assouad dimension
  is again $\frac{\log4}{\log 3}$, see \cite[Theorem~2.2]{MR3267023} and \cite[Theorem~A.5(3)]{MR1608518}. Hence $\text{co\,dim}_{A}(\partial\Omega)=2-\frac{\log4}{\log 3}$.

  Moreover, by \cite[Theorem 1.1]{MR2269586} the volume of the inner tubular neighbourhood of $\Omega$ is described by the formula
$$
\left|\Omega_{r}\right|=G_{1}(r)r^{2-\frac{\log4}{\log3}}+G_{2}(r)r^{2},
$$
where $G_{1}$ and $G_{2}$ are continuous, periodic functions (in consequence bounded). Hence, for $r<1$ we have $\left|\Omega_{r}\right|=O\left(r^{2-\frac{\log4}{\log3}}\right)$. Since in addition $\Omega$ is $\kappa$ - plump, by Theorem \ref{tw1} we obtain that if $p=1$, then $C_{c}^{\infty}(\Omega)$ is dense in $W^{s,p}(\Omega)$ if $s< 2-\frac{\log 4}{\log 3}$ and is not dense if $s> 2-\frac{\log 4}{\log 3}$. Moreover, if $p>1$, then the density result holds if and only if $sp\leq 2-\frac{\log4}{\log3}$. We do not know what is happening in the remaining case $p=1$ and $s=2-\frac{\log4}{\log3}$.
\end{example}

\section{The space $W_{0}^{s,p}(\Omega)$}\label{sec:W0}
Based on our previous results, we are able to describe explicitly the space $W_{0}^{s,p}(\Omega)$ in some particular cases. Namely, we can describe this space for $\Omega, s$ and $p$ satisfying the following weak fractional Hardy inequality.
	\begin{defn}\label{defn5}
	We say that $\Omega$ admits a weak $(s,p)$--fractional Hardy inequality, if there exists a constant $c=c(d,s,p,\Omega)$ such that for every $f\in C_{c}^{\infty}(\Omega)$ it holds
\[
\int_{\Omega}\frac{|f(x)|^{p}}{d_{\Omega}(x)^{sp}}\,dx\leq c \| f\|_{W^{s,p}(\Omega)}^p.
\]
In the case when the norm $\| f\|_{W^{s,p}(\Omega)}$ above can be replaced by the seminorm
$[ f]_{W^{s,p}(\Omega)}$, we say that $\Omega$ admits an $(s,p)$--fractional Hardy inequality.
	\end{defn}
\begin{thm}\label{tw3}
    Suppose that $\Omega$ admits a weak $(s,p)$-fractional Hardy inequality. Then
    $$
    W_{0}^{s,p}(\Omega)=\left\{f\in W^{s,p}(\Omega):\int_{\Omega}\frac{|f(x)|^{p}}{d_{\Omega}(x)^{sp}}\,dx<\infty\right\}.
    $$
\end{thm}
\begin{proof}
By Lemma \ref{lem1}, if $\int_{\Omega}\frac{|f(x)|^{p}}{d_{\Omega}(x)^{sp}}\,dx<\infty$, then $f\in W_{0}^{s,p}(\Omega)$, because in this case
$$
n^{sp}\int_{\Omega_{\frac{3}{n}}}|f(x)|^{p}\,dx\leq 3^{sp}\int_{\Omega_{\frac{3}{n}}}\frac{|f(x)|^{p}}{d_{\Omega}(x)^{sp}}\,dx\longrightarrow 0,
$$
when $n\longrightarrow\infty$. In fact, for that part we do not need the assumption about Hardy inequality.

 Suppose that $\Omega$ admits a weak $(s,p)$--Hardy inequality and $f\in W_{0}^{s,p}(\Omega)$. Let $f_{n}$ be a sequence of smooth and compactly supported functions convergent to $f$ in $W^{s,p}(\Omega).$ In particular, $f_{n}\longrightarrow f$ in $L^{p}(\Omega)$, so there exists a subsequence $f_{n_{k}}$ convergent to $f$ almost everywhere. We have by Fatou lemma
\begin{align*}
\int_{\Omega}\frac{|f(x)|^{p}}{d_{\Omega}(x)^{sp}}\,dx&=\int_{\Omega}\lim_{k\rightarrow\infty}\frac{|f_{n_{k}}(x)|^{p}}{d_{\Omega}(x)^{sp}}\,dx\\
&\leq\liminf_{k\rightarrow\infty}\int_{\Omega}\frac{|f_{n_{k}}(x)|^{p}}{d_{\Omega}(x)^{sp}}\,dx\\
&\leq c\liminf_{k\rightarrow\infty}\|f_{n_{k}}\|^{p}_{W^{s,p}(\Omega)}\\
&=c\|f\|^{p}_{W^{s,p}(\Omega)} < \infty. \qedhere
\end{align*}
\end{proof}

\begin{proof}[Proof of Theorem~\ref{twW0}]
  From part (F) of Theorem~\ref{thm.dv} with $\eta=sp$, $\varphi(t)=t^{sp}$, $\Omega$ admits an $(s,p)$-fractional Hardy inequality and also a weak $(s,p)$-fractional Hardy inequality. Thus the result follows from Theorem~\ref{tw3}.
\end{proof}

\begin{proof}[Proof of Theorem~\ref{weakH}]
  From part (T') of Theorem~\ref{thm.dv}, inequality \eqref{eq:weakH} holds for all functions $f$ for which the left hand side of \eqref{eq:weakH} is finite.
  Thus by Theorem~\ref{twW0}, it holds for all functions $f\in W_{0}^{s,p}(\Omega)$.
 However, by part (I) of Theorem~\ref{tw1},  $W_{0}^{s,p}(\Omega) =  W^{s,p}(\Omega)$
  and the result follows.
\end{proof}

\section{Appendix}
We recall from \cite[Section 3]{MR3165234} the notion of a global weak lower (or upper) scaling condition ($\WLSC$ or $\WUSC$ for short). As in \cite{MR3237044}, we will use a different, but equivalent formulation. We note that in our setting the middle parameter in $\WLSC$ or $\WUSC$ is always zero and thus we could omit it, however we prefer to keep the notation consistent with \cite{MR3165234, MR3237044}.

\begin{defn}\label{d.WUSC}
Let $\eta\in \R$ and  $H\in (0,1]$.
We say that a function $\phi:(0,\infty)\to(0,\infty)$ satisfies $\WLSC(\eta, 0, H)$
(resp.,  $\WUSC(\eta, 0, H^{-1})$)
and write $\phi \in \WLSC(\eta, 0, H)$ ($\phi \in \WUSC(\eta, 0, H^{-1})$), if
\begin{align}
\phi(st) \geq H t^\eta \phi(s), \quad s>0\,, \label{eq:phi-assumption}
\end{align}
for every $t\geq 1$ (resp., for every $t\in (0,1]$).
\end{defn}

We begin with the following observation:
\begin{equation}\label{eq.observation}
  \text{If $\Omega\subset \R^d$ is a nonempty open bounded set, then $\udimA (\partial \Omega) \geq d-1$.}
\end{equation}
For the proof, we will provide the following argument by the user \rm{rpotrie} from  \cite{MathOverflow}.
Since $\partial \Omega$ disconnects $\R^d$, its topological dimension has to be at least $d-1$, see \cite[Theorem IV.4]{MR0006493}. But the topogical dimension does not exceed Hausdorff dimension \cite[page 107]{MR0006493}, and the latter in turn does not exceed the upper Assouad dimension  \cite[Theorem A.5(10)]{MR1608518}, consequently \eqref{eq.observation} holds.

\begin{proof}[Proof of case (T') in Theorem~\ref{thm.dv}]
  It seems possible to adapt the original proof for this case, however, since the proof was quite involved and technical, we prefer to choose another strategy. Namely, we will reduce (T') to the case (T).
  Let us assume that the general assumptions of Theorem~\ref{thm.dv} and the assumptions in (T') hold.
  
  Let us fix $x_0 \in \Omega$ and put $M=\diam \Omega$.
  We consider an open set $\Omega_1 = \R^d \setminus \overline B(x_0, 2M)$. Let $G = \Omega\cup \Omega_1$.
  Observe that $\dist(\Omega, \Omega_1) \geq M$, hence $\partial G = \partial \Omega \cup \partial \Omega_1$.
  Therefore,
  \[
  \udimA (\partial{G}) = \max\{ \udimA (\partial \Omega),\, \udimA (\partial \Omega_1) \}
  = \max\{ \udimA (\partial \Omega), d-1 \} =  \udimA (\partial \Omega),
  \]
  by \cite[Theorem A.5(3)]{MR1608518} and \eqref{eq.observation}.

  We may also need to redefine the function $\phi$. To this end, put $\eta_0 = \eta$ if $\eta>0$, while in the case when $\eta\leq 0$, we choose $\eta_0>0$ such that
  \[
  \eta_0 + \overline{\mathrm{dim}}_A(\partial \Omega) - d<0.
  \]
  We note that this is possible, because $\kappa$-plumpness of $\Omega$ implies that $\partial \Omega$
  is porous, and that in turn by \cite[Theorem 5.2]{MR1608518} implies that $\overline{\mathrm{dim}}_A(\partial \Omega) < d$.
  We define
  \[
  \psi(x) = \begin{cases}
    \phi(x),&\text{when $x \in (0, M]$};\\
      \phi(T)(\frac{x}{T})^{\eta_0},&\text{when $x \in (M, \infty)$}.
  \end{cases}
  \]
  We claim that such a~function $\psi$ satisfies the condition $\WUSC(\eta_0,0,H^{-1})$. We omit a straightforward check of \eqref{eq:phi-assumption} in three possible cases, when the two numbers $st \leq s$ in that equation  lie in either $(0, M]$ or $(M, \infty)$.

    We apply the case (T) of the Theorem~\ref{thm.dv} (proved in \cite{MR3237044}) to the open set $G$, the number $\eta_0$ and the function $\psi \in \WUSC(\eta_0,0,H^{-1})$. It follows that there exist constants $c$ and $R$ such that
    \begin{equation}\label{eq.omegaHardy}
 \int_G \frac{|u(x)|^p}{\psi(d_{G}(x))}\,dx \leq
c \int_G\!\int_{G\cap B(x,Rd_{G}(x))}
\frac{|u(x)-u(y)|^p}{\psi(d_{G}(x))d_{G}(x)^d}\,dy\,dx
\end{equation}
holds for all measurable functions $u:G\to \R$ for which the left hand side is finite.

Let us consider an arbitrary measurable functions $u:\Omega\to \R$ for which $\int_\Omega \frac{|u(x)|^p}{\phi(d_{G}(x))}\,dx < \infty$, and extend it by zero on $\Omega_1$ to obtain a function defined on the whole set $G$. Inequality \eqref{eq.omegaHardy} for this function $u$ has the following form,
\begin{align*}
  \int_\Omega \frac{|u(x)|^p}{\phi(d_{G}(x))}\,dx &\leq
c \int_\Omega\!\int_{\Omega\cap B(x,Rd_{G}(x))} 
\frac{|u(x)-u(y)|^p}{\phi(d_{G}(x))d_{G}(x)^d}\,dy\,dx \\
&+c\int_{\Omega_1}\!\int_{\Omega\cap B(x,Rd_{G}(x))}
\frac{|u(y)|^p}{\psi(d_{G}(x))d_{G}(x)^d}\,dy\,dx \\
&+c \int_{\Omega}\!\int_{\Omega_1\cap B(x,Rd_{G}(x)))}
\frac{|u(x)|^p}{\psi(d_{G}(x))d_{G}(x)^d}\,dy\,dx\\
&=:c(I_1+I_2+I_3).
\end{align*}

In the integral $I_2$, when $x\in \Omega_1$ and $y\in \Omega\cap B(x,Rd_{G}(x))$, then
$M\leq |x-y| \leq Rd_{G}(x)$ and therefore $d_{G}(x) \geq M/R$. Consequently,
\begin{equation}\label{eq.I2}
I_2\leq \|u\|_{L^p(\Omega)}^p \int_{\{x\in \Omega_1:d_{G}(x) \geq M/R \}} \frac{dx}{\psi(d_{G}(x))d_{G}(x)^d}.
\end{equation}

From the definition of the function $\psi$ and the fact that $\psi\in \WUSC(\eta_0,0,H^{-1})$
it follows that there exists a constant $c(M/R, H, \eta_0)$ such that
\begin{equation}\label{eq.psi_asymp}
\psi(z) \geq c(M/R, H, \eta_0) z^{\eta_0}, \quad \text{for $z\geq M/R$.}
\end{equation}

Therefore the integral in \eqref{eq.I2} is convergent and so $I_2 \leq c'  \|u\|_{L^p(\Omega)}^p$.

For the integral $I_3$ we observe that when $x\in \Omega$ and $y\in \Omega_1\cap B(x,Rd_{G}(x))$,
then $d_G(x)=d_\Omega(x)$ and $M\leq |y-x| \leq Rd_{G}(x)$, so $d_{G}(x) \geq M/R$. Therefore by \eqref{eq.psi_asymp}
the function $\psi(d_{G}(x))^{-1} d_{G}(x)^{-d}$ is bounded from above.
Furthermore, since $|y-x_0| \leq M+|y-x| \leq M+Rd_{G}(x) \leq M(1+R)$, the following inclusion
$\Omega_1\cap B(x,Rd_{G}(x)) \subset B(x_0, M(1+R))$  holds for all $x\in \Omega$.
Thus also in this case $I_3 \leq c'  \|u\|_{L^p(\Omega)}^p$.

Consequently, $I_1$ is equal to the first term on the right side of \eqref{eq:fhi},
while $I_2$ and $I_3$ are bounded by the second term.
\end{proof}


\begin{thebibliography}{10}
\bibitem{MR3989177}
{\sc Baalal, A., and Berghout, M.}
\newblock Density properties for fractional {S}obolev spaces with variable
  exponents.
\newblock {\em Ann. Funct. Anal. 10}, 3 (2019), 308--324.

\bibitem{MR3165234}
{\sc Bogdan, K., Grzywny, T., and Ryznar, M.}
\newblock Density and tails of unimodal convolution semigroups.
\newblock {\em J. Funct. Anal. 266}, 6 (2014), 3543--3571.

\bibitem{MR3309383}
{\sc Botelho, F.}
\newblock {\em Functional analysis and applied optimization in {B}anach
  spaces}.
\newblock Springer, Cham, 2014.
\newblock Applications to non-convex variational models, With contributions by
  Anderson Ferreira and Alexandre Molter.

\bibitem{MR4225499}
{\sc Brasco, L., G\'{o}mez-Castro, D., and V\'{a}zquez, J.~L.}
\newblock Characterisation of homogeneous fractional {S}obolev spaces.
\newblock {\em Calc. Var. Partial Differential Equations 60}, 2 (2021), Paper
  No. 60, 40.

\bibitem{MR3987216}
{\sc Brasco, L., and Salort, A.}
\newblock A note on homogeneous {S}obolev spaces of fractional order.
\newblock {\em Ann. Mat. Pura Appl. (4) 198}, 4 (2019), 1295--1330.

\bibitem{MR1792288}
{\sc Caetano, A.~M.}
\newblock Approximation by functions of compact support in
  {B}esov-{T}riebel-{L}izorkin spaces on irregular domains.
\newblock {\em Studia Math. 142}, 1 (2000), 47--63.

\bibitem{MR1993864}
{\sc Chen, Z.-Q., and Song, R.}
\newblock Hardy inequality for censored stable processes.
\newblock {\em Tohoku Math. J. (2) 55}, 3 (2003), 439--450.

\bibitem{MR3420496}
{\sc Dipierro, S., and Valdinoci, E.}
\newblock A density property for fractional weighted {S}obolev spaces.
\newblock {\em Atti Accad. Naz. Lincei Rend. Lincei Mat. Appl. 26}, 4 (2015),
  397--422.

\bibitem{MR2085428}
{\sc Dyda, B.}
\newblock A fractional order {H}ardy inequality.
\newblock {\em Illinois J. Math. 48}, 2 (2004), 575--588.

\bibitem{MR4190640}
{\sc Dyda, B., and Kijaczko, M.}
\newblock On density of smooth functions in weighted fractional {S}obolev
  spaces.
\newblock {\em Nonlinear Anal. 205\/} (2021), 112231, 10.

\bibitem{MR3237044}
{\sc Dyda, B., and V\"{a}h\"{a}kangas, A.~V.}
\newblock A framework for fractional {H}ardy inequalities.
\newblock {\em Ann. Acad. Sci. Fenn. Math. 39}, 2 (2014), 675--689.

\bibitem{MR3310082}
{\sc Fiscella, A., Servadei, R., and Valdinoci, E.}
\newblock Density properties for fractional {S}obolev spaces.
\newblock {\em Ann. Acad. Sci. Fenn. Math. 40}, 1 (2015), 235--253.

\bibitem{MR3267023}
{\sc Fraser, J.~M.}
\newblock Assouad type dimensions and homogeneity of fractals.
\newblock {\em Trans. Amer. Math. Soc. 366}, 12 (2014), 6687--6733.

\bibitem{MR775683}
{\sc Grisvard, P.}
\newblock {\em Elliptic problems in nonsmooth domains}, vol.~24 of {\em
  Monographs and Studies in Mathematics}.
\newblock Pitman (Advanced Publishing Program), Boston, MA, 1985.

\bibitem{MR4144553}
{\sc Henderson, A.~M.}
\newblock {\em Fractal {Z}eta {F}unctions in {M}etric {S}paces}.
\newblock ProQuest LLC, Ann Arbor, MI, 2020.
\newblock Thesis (Ph.D.)--University of California, Riverside.

\bibitem{MR0006493}
{\sc Hurewicz, W., and Wallman, H.}
\newblock {\em Dimension {T}heory}.
\newblock Princeton Mathematical Series, vol. 4. Princeton University Press,
  Princeton, N. J., 1941.

\bibitem{MR3277052}
{\sc Ihnatsyeva, L., Lehrb\"{a}ck, J., Tuominen, H., and V\"{a}h\"{a}kangas,
  A.~V.}
\newblock Fractional {H}ardy inequalities and visibility of the boundary.
\newblock {\em Studia Math. 224}, 1 (2014), 47--80.

\bibitem{MR3205534}
{\sc K\"{a}enm\"{a}ki, A., Lehrb\"{a}ck, J., and Vuorinen, M.}
\newblock Dimensions, {W}hitney covers, and tubular neighborhoods.
\newblock {\em Indiana Univ. Math. J. 62}, 6 (2013), 1861--1889.

\bibitem{MR1470421}
{\sc Kinnunen, J., and Martio, O.}
\newblock Hardy's inequalities for {S}obolev functions.
\newblock {\em Math. Res. Lett. 4}, 4 (1997), 489--500.

\bibitem{MR664599}
{\sc Kufner, A.}
\newblock {\em Weighted {S}obolev spaces}, vol.~31 of {\em Teubner-Texte zur
  Mathematik [Teubner Texts in Mathematics]}.
\newblock BSB B. G. Teubner Verlagsgesellschaft, Leipzig, 1980.
\newblock With German, French and Russian summaries.

\bibitem{MR2269586}
{\sc Lapidus, M.~L., and Pearse, E. P.~J.}
\newblock A tube formula for the {K}och snowflake curve, with applications to
  complex dimensions.
\newblock {\em J. London Math. Soc. (2) 74}, 2 (2006), 397--414.

\bibitem{MR1608518}
{\sc Luukkainen, J.}
\newblock Assouad dimension: antifractal metrization, porous sets, and
  homogeneous measures.
\newblock {\em J. Korean Math. Soc. 35}, 1 (1998), 23--76.

\bibitem{MR1810244}
{\sc Mamedov, F.~I.}
\newblock On the multidimensional weighted {H}ardy inequalities of fractional
  order.
\newblock {\em Proc. Inst. Math. Mech. Acad. Sci. Azerb. 10\/} (1999),
  102--114, 275.

\bibitem{MathOverflow}
{\sc MathOverflow}.
\newblock {H}ausdorff dimension of the boundary of an open set in the
  {E}uclidean space -- lower bound;
  \url{https://mathoverflow.net/questions/40593/hausdorff-dimension-of-the-boundary-of-an-open-set-in-the-euclidean-space-lowe};
  the answer of the user \rm{rpotrie}.

\bibitem{MR1742312}
{\sc McLean, W.}
\newblock {\em Strongly elliptic systems and boundary integral equations}.
\newblock Cambridge University Press, Cambridge, 2000.

\bibitem{MR3024598}
{\sc Triebel, H.}
\newblock {\em Theory of function spaces}.
\newblock Modern Birkh\"{a}user Classics. Birkh\"{a}user/Springer Basel AG,
  Basel, 2010.
\newblock Reprint of 1983 edition.


\bibitem{10.2748/tmj/1178228081}
{\sc Väisälä, J.}
\newblock {Uniform domains}.
\newblock {\em Tohoku Mathematical Journal 40}, 1 (1988), 101 -- 118.

\end{thebibliography}
\end{document}